\documentclass{article}

\usepackage{longtable}
\usepackage{charter}
\usepackage{geometry}
\usepackage{graphicx}
\usepackage{amsmath,amssymb}
\usepackage{xcolor}
\usepackage{youngtab}
\usepackage{pifont}

\usepackage[colorlinks=true,bookmarks=false]{hyperref}
\definecolor{myurlcolor}{rgb}{0.6,0,0}
\definecolor{mycitecolor}{rgb}{0,0,0.8}
\definecolor{myrefcolor}{rgb}{0,0,0.8}
\hypersetup{linkcolor=myrefcolor,citecolor=mycitecolor,urlcolor=myurlcolor}

\definecolor{lblue}{HTML}{cce6ff}
\definecolor{lyellow}{HTML}{e0ba4f}
\definecolor{tan}{HTML}{fcba03}
\definecolor{purple}{HTML}{bd358f}
\definecolor{lightgreen}{HTML}{7fde31}
\definecolor{darkorange}{HTML}{b55a2a}
\definecolor{darkgreen}{HTML}{4e8f18}
\definecolor{lightblue}{rgb}{0.9,0.95,1}
\definecolor{darkblue}{rgb}{0,0,100}
\definecolor{brickred}{rgb}{0.7,0,0}

\newcommand{\define}[1]{{\bf \boldmath{#1}}}

\newcommand\R{{\mathbb R}}
\newcommand\C{{\mathbb C}}

\newcommand{\U}{{\mathrm U }}
\newcommand{\SU}{{\mathrm{SU}}}
\newcommand{\GL}{{\mathrm{GL}}}
\newcommand{\SL}{{\mathrm{SL}}}

\newcommand{\End}{{\mathrm{End}}}

\newcommand{\maps}{\colon}
\newcommand{\RR}{\mathsf{R}}
\newcommand{\Rep}{\mathsf{Rep}}
\newcommand{\Schur}{\mathsf{Schur}}

\renewcommand{\texttt}[1]{%
  \begingroup
  \ttfamily
  \begingroup\lccode`~=`/\lowercase{\endgroup\def~}{/\discretionary{}{}{}}%
  \begingroup\lccode`~=`[\lowercase{\endgroup\def~}{[\discretionary{}{}{}}%
  \begingroup\lccode`~=`.\lowercase{\endgroup\def~}{.\discretionary{}{}{}}%
  \catcode`/=\active\catcode`[=\active\catcode`.=\active
  \scantokens{#1\noexpand}%
  \endgroup
}

\usepackage{tikz}
\usepackage{tikz-cd}
\usetikzlibrary{matrix,arrows,decorations.pathmorphing,positioning}
\usetikzlibrary{intersections,decorations.markings}
\usetikzlibrary{arrows,positioning,fit,matrix,shapes.geometric,external}
\usetikzlibrary{backgrounds,circuits,circuits.ee.IEC,shapes,fit,matrix}

\tikzstyle{inarrow}=[->, >=stealth, shorten >=.03cm,line width=1.5]

\makeatletter
\newcommand{\xRightarrow}[2][]{\ext@arrow 0359\Rightarrowfill@{#1}{#2}}
\makeatother

\begin{document}

\thispagestyle{empty}

  \begin{center}
    {\Huge\textbf{Young Diagrams and \\
    Classical Groups}}  
 \\[4em]
  {John Baez}
  \\
  {Department of Mathematics}
  \\
  {University of California, Riverside}
  \\[0.5em]{\today}
  \\[0.5em]{{based on ``week157'' of \emph{\href{http://math.ucr.edu/home/baez/TWF.html'}{This Week's Finds}}}}
  \\[2em]
  \includegraphics[scale = 0.7]{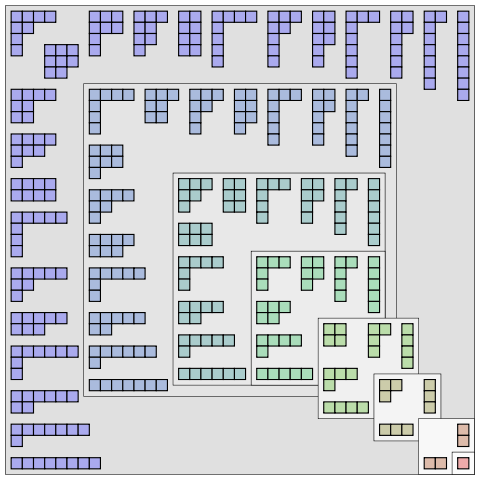} \\
  \tiny{\href{https://commons.wikimedia.org/wiki/File:Ferrer_partitioning_diagrams.svg}{image by R.\ A.\ Nonenmacher, CC BY-SA 4.0}}
  \end{center}
  
\newpage
\setcounter{page}{1}

Mathematics and physics rely a lot on \emph{symmetry} to simplify problems, and there
are two kinds of diagrams that show up a lot in this context: Dynkin
diagrams and Young diagrams. Dynkin diagrams first show up when you study shapes with lots of
reflection symmetries, like crystals and Platonic solids. They wind up
being good for all sorts of other stuff, like classifying simple Lie groups and their representations.  But what about Young diagrams? These are also important for studying group representations, but for a more limited class of groups: the ``classical'' groups.    Representations of classical groups are used a lot in quantum physics, from particle physics through nuclear physics and atomic physics up to chemistry.    So Young diagrams are not only beautiful, they're practical.

My goal is to explain how Young diagrams are used to classify representations of classical groups.  I won't prove much, just sketch the ideas.  First I'll explain classical groups and group representations.   But even before that, I should say what's a Young diagram. 

\subsubsection*{Young diagrams}

Here is an example of a Young diagram:
\[    \yng(6,5,5,2,1) \]
All the information here is captured by the number of boxes in each row:
\[   6 \ge 5 \ge 5 \ge 2 \ge 1 \]
So, we can define a \define{Young diagram} to be a finite sequence of natural numbers \(n_1 \ge n_2 \ge \cdots \ge n_k > 0$.    We say \(k\) is the number of \define{rows} and  \(n_1\) is the number of \define{columns}.  We say \(n_i\) is the number of \define{boxes} in the \(i$th column, and \(n = \sum_i n_i\) is the total number of boxes.

Young diagrams with \(n\) boxes classify partitions of an \(n$-element set, up to isomorphism.  For example, this partition:
\[
\scalebox{0.6}{
\begin{tikzpicture}
\path
  coordinate (aux0) at (0,1.5)
  coordinate (aux1) at (0,3.5)
  coordinate (aux2) at (10,3.5)
  coordinate (aux3) at (9,6)
  coordinate (aux4) at (4,0)
  coordinate (aux5) at (7,0)
  coordinate (aux6) at (2,6)
  coordinate (aux7) at (5,6)
  coordinate (esp1) at (0.2,2.5)
  coordinate (esp2) at (1.5,1.5)
  coordinate (esp3) at (3,0.1)
  coordinate (esp4) at (5.5,1.1)
  coordinate (esp5) at (8,0.5)
  coordinate (esp6) at (8.75,2)
  coordinate (esp7) at (9.7,3)
  coordinate (esp8) at (6.5,4.5)
  coordinate (esp9) at (3.8,5.8)
  coordinate (esp10) at (1.5,4)
  ;
\draw[line width=0.8pt]
  (esp1) to[out=-90,in=170]
  (esp2) to[out=-10,in=170]
  (esp3) to[out=-10,in=180]
  (esp4) to[out=0,in=180]
  (esp5) to[out=10,in=-150]
  (esp6) to[out=20,in=-90]
  (esp7) to[out=90,in=-60]
  (esp8) to[out=120,in=0]
  (esp9) to[out=180,in=0]
  (esp10) to[out=180,in=90]
  cycle;    
\clip
  (esp1) to[out=-90,in=170]
  (esp2) to[out=-10,in=170]
  (esp3) to[out=-10,in=180]
  (esp4) to[out=0,in=180]
  (esp5) to[out=10,in=-150]
  (esp6) to[out=20,in=-90]
  (esp7) to[out=90,in=-60]
  (esp8) to[out=120,in=0]
  (esp9) to[out=180,in=0]
  (esp10) to[out=180,in=90]
  cycle;    
\filldraw[fill=darkorange!40]
  (aux4) to[bend right=10]
  (aux6) --
  (aux7) to[bend left=10]
  (aux5) -- cycle;
\filldraw[fill=brown!60]
  (aux5) to[bend right=10]
  (aux7) --
  (10,6) --
  (10,0) -- cycle;
\filldraw[fill=tan!30]
  (aux0) -- 
  (aux1) to[bend right=10]
  (aux3) --
  (10,6) -- 
  (aux2) to[bend left=10] cycle;
\filldraw[fill=lyellow!50]
  (0,0) -- 
  (aux4) to[bend right=10]
  (aux6) --
  (0,6) -- 
  (0,0) -- cycle;
\filldraw[fill=orange!40]
  (0,6) -- 
  (aux1) to[bend right=10]
  (aux3) --
  (0,6) -- cycle;
\node at (3,5) {$\bullet$};  
\node at (5,5.2) {$\bullet$};  
\node at (2,2) {$\bullet$};  
\node at (1,3) {$\bullet$};  
\node at (3,2.8) {$\bullet$};  
\node at (5,4) {$\bullet$};  
\node at (4,3) {$\bullet$};  
\node at (6.4,3.4) {$\bullet$};
\node at (4.4,1.3) {$\bullet$};  
\node at (6,1.7) {$\bullet$};  
\node at (7.5,1) {$\bullet$};  
\node at (8.2,2) {$\bullet$};  
\end{tikzpicture}
}
\]
gives this Young diagram, whose rows list how many points are in each part:
\[  \yng(3,3,2,2,2)  \]
But the Young diagram does not record which point of our set lies in which part, so
Young diagrams classify partitions only ``up to isomorphism''.

Young diagrams with \(n\) boxes also classify permutations of an \(n$-element set
up to isomorphism.  For example this permutation: 
\[ 
\color{darkgreen}
\scalebox{0.6}{
\begin{tikzpicture}
\path
  coordinate (aux0) at (0,1.5)
  coordinate (aux1) at (0,3.5)
  coordinate (aux2) at (10,3.5)
  coordinate (aux3) at (9,6)
  coordinate (aux4) at (4,0)
  coordinate (aux5) at (7,0)
  coordinate (aux6) at (2,6)
  coordinate (aux7) at (5,6)
  coordinate (esp1) at (0.2,2.5)
  coordinate (esp2) at (1.5,1.5)
  coordinate (esp3) at (3,0.1)
  coordinate (esp4) at (5.5,1.1)
  coordinate (esp5) at (8,0.5)
  coordinate (esp6) at (8.75,2)
  coordinate (esp7) at (9.7,3)
  coordinate (esp8) at (6.5,4.5)
  coordinate (esp9) at (3.8,5.8)
  coordinate (esp10) at (1.5,4)
  ;
\draw[line width=1.2pt]
  (esp1) to[out=-90,in=170]
  (esp2) to[out=-10,in=170]
  (esp3) to[out=-10,in=180]
  (esp4) to[out=0,in=180]
  (esp5) to[out=10,in=-150]
  (esp6) to[out=20,in=-90]
  (esp7) to[out=90,in=-60]
  (esp8) to[out=120,in=0]
  (esp9) to[out=180,in=0]
  (esp10) to[out=180,in=90]
  cycle;    
\clip
  (esp1) to[out=-90,in=170]
  (esp2) to[out=-10,in=170]
  (esp3) to[out=-10,in=180]
  (esp4) to[out=0,in=180]
  (esp5) to[out=10,in=-150]
  (esp6) to[out=20,in=-90]
  (esp7) to[out=90,in=-60]
  (esp8) to[out=120,in=0]
  (esp9) to[out=180,in=0]
  (esp10) to[out=180,in=90]
  cycle;    
\node at (3,5) (A) {$\bullet$};  
\node at (5,5.2) (B) {$\bullet$};  
\node at (2,2) (C) {$\bullet$};  
\node at (1,3) (D) {$\bullet$}; 
\node at (3,2.8) (E) {$\bullet$};  
 
\node at (5,4) (F) {$\bullet$};  
\node at (4,3) (G) {$\bullet$};  
\node at (6.4,3.4) (H) {$\bullet$};
\node at (4.4,1.3) (I) {$\bullet$};  
\node at (6,1.7) (J) {$\bullet$};  
\node at (7.5,1) (K) {$\bullet$};  
\node at (8.2,2) (L) {$\bullet$};  

\draw [style=inarrow, bend left=20, looseness=1.00] (A) to (B);
\draw [style=inarrow, bend left=20, looseness=1.00] (B) to (A);
\draw [style=inarrow, bend left=30, looseness=1.00] (C) to (D);
\draw [style=inarrow, bend left=30, looseness=1.00] (D) to (E);
\draw [style=inarrow, bend left=30, looseness=1.00] (E) to (C);
\draw [style=inarrow, bend right=30, looseness=1.00] (F) to (G);
\draw [style=inarrow, bend right=30, looseness=1.00] (G) to (H);
\draw [style=inarrow, bend right=30, looseness=1.00] (H) to (F);
\draw [style=inarrow, bend right=30, looseness=1.00] (I) to (J);
\draw [style=inarrow, bend right=50, looseness=1.00] (J) to (I);
\draw [style=inarrow, bend right=30, looseness=1.00] (K) to (L);
\draw [style=inarrow, bend right=30, looseness=1.00] (L) to (K);
\end{tikzpicture}
}\]
gives the same Young diagram we have just seen.  But any isomorphic permutation
would give the same Young diagram.

What's an ``isomorphic permutation'', exactly?   Let's look at an example.   Permutations of the set \(\{1,\dots, n\}\) form the \define{symmetric group} \(S_n$.   Say we have any permutation \(g \in S_n\), like this: 
\[
  \begin{aligned}
    1 &\to 2
  \\2 &\to 4 
  \\3 &\to 3
  \\4 &\to 1
  \\5 &\to 6
  \\6 &\to 5
  \\7 &\to 7  
  \end{aligned}
\] 
Note that 1 gets mapped to 2, which gets mapped to 4, which gets
mapped back to 1 again. Similarly, 5 gets mapped to 6, which gets mapped
back to 5. The number 3 gets mapped to itself right away, as does 7. No
matter where we start, we always cycle back eventually. So our
permutation consists of a bunch of cycles: 
\[(1,2,4) (5,6) (3) (7)\]
and this ``cycle decomposition'' completely describes the
permutation. To simplify life, we always write down these cycles in
order of decreasing length. We also write the lowest number in each
cycle first.

Now suppose we conjugate our permutation \(g\) by some other
permutation, say \(h\). This gives the permutation \(hgh^{-1}\). How
does the cycle decomposition of this compare with that of \(g\)? It
looks very similar!  For example, it might look like this:
\[(2,7,6) (1,3) (4) (5)\] 
There are the same number of cycles, each the
same length as before. The only thing that changes are the numbers in
each cycle. These get switched around by means of the permutation \(h\).

In short, when we conjugate a permutation, all that remains unchanged is
the picture we get by writing down its cycle decomposition and blotting
out the specific numbers in each cycle, like this:
\[(\square,\square,\square) (\square,\square) (\square) (\square)\]
If we write each cycle as a row of boxes, we get a Young diagram:
\[
 \yng(3,2,1,1)
\] 

\subsubsection*{Classical groups, and a classical monoid}

Now, what are the classical groups? 
As with composers of music, there's no precise list of groups that count
as ``classical''. But in general, a classical group should consist of
linear transformations that preserve some nice geometrical structure on
a vector space.   Some good examples are:
\begin{itemize}
\item
The \define{general linear group} \(\GL(N,\C)\), consisting of all invertible
linear transformations of \(\C^N$, or in other words, all \(N \times N\) complex
matrices with nonzero determinant.
\item 
The \define{special linear group} \(\SL(N,\C)\), consisting of all
linear transformations of \(\C^N\) with determinant \(1\).   
\item 
The \define{unitary group} \(\U(N)\), consisting of all unitary linear
transformations of \(\C^N\).
\item 
The \define{special unitary group} \(\SU(N)\), consisting of all unitary linear
transformations of \(\C^N\) with determinant \(1\).
\end{itemize}
These are the Bach, Haydn, Mozart and Beethoven of classical groups.  Representations of all 
four can be classified with the help of Young diagrams.  

We may also consider this an honorary classical group, even though it's defined in terms of a \emph{set} rather than a \emph{vector space}:
\begin{itemize}
\item 
The symmetric group \(S_n\), consisting of all permutations of the set \( \{1,\dots,n\} \).
\end{itemize}
Representations of this group are also classified using Young diagrams---and as we'll
see, \(S_n\) plays a starring role in the whole story.

There's another key actor whose representations are classified by Young diagrams.  It deserves to be called a ``classical monoid'':
\begin{itemize}
\item 
The \define{full linear monoid} $\End(\C^n)$, consisting of all linear transformations of \(\C^N$, or in other words, all \(N \times N\) matrices.
\end{itemize}
A \define{monoid} is a set with an associative multiplication and identity,
but not necessarily inverses.  Here I am making \(\End(\C^n)\) into a monoid where the multiplication is composition of transformations---or in low-brow terms, matrix multiplication.
This monoid is so classical that people don't even call it that!  Perhaps the common prejudice in favor of groups and against other monoids is to blame.  As we'll see, the full linear monoid is a bit like the composer Palestrina, who is not considered a classical composer, yet who set the stage for the music we call classical.

\subsubsection*{\bf Representations}

Groups feel sad unless they are acting as symmetries of something.   Monoids feel the same way---or even worse, because they're less loved than groups.   This why we should study representations of groups and monoids.   A \define{homomorphism} of monoids, say \(\rho \maps M \to N$, is a function with
\[       \rho(mm') = \rho(m) \rho(m') \textrm{ for all } m, m' \in M \textrm{ and } \rho(1) = 1 .\]
A \define{representation} of a monoid \(M\) on a vector space \(V\) is a homomorphism
\[       \rho \maps M \to \End(V) \]
where \(\End(V)\) consists of all linear transformations of \(V\), made into a monoid using
composition.  A representation lets us take an element \(m \in M\) and make it act on a vector \(v \in V\) to get a new vector \(\rho(m)v$, in such a way that
\[      \rho(mm') v = \rho(m) \rho(m') v \textrm{ and } \rho(1) v = v .\]
So now our monoid is doing something, not just sitting there moping!

But a representation is still lonely in isolation.  To solve this problem we define morphisms between representations of given monoid, getting an entire \emph{category} of representations.   
Given two representations \(\rho \maps M \to \End(V), \sigma \maps
M \to \End(W)\), a \define{morphism} from the first to the second is
a linear map \(f \maps V \to W\) such that
\[     f(\rho(m) v) = \sigma(m) f(v)   \]
for all \(v \in V\).   That is: acting and then mapping is the same as mapping and then
acting.   Thanks to how \(f\) slips from outside to inside in this equation,
morphisms of representations are also called \define{intertwining operators}.   

An \define{isomorphism} is just a morphism with an inverse, and an isomorphism of 
representations is also commonly called an \define{equivalence}.   We won't do much
with categories here except for classifying representations ``up to isomorphism'': 
when we do that, we don't distinguish between isomorphic representations.   But studying
the whole category of representations of a monoid, all at once, is a good way to get 
deeper insights in representation theory.

The simplest representations are those on finite-dimensional vector spaces---so henceforth:
\begin{center}
{\color{brickred} We assume all vector spaces under discussion 
are finite-dimensional, without even mentioning it!}
\end{center}
And instead of trying to study \emph{all} finite-dimensional representations, I will
focus on the ``irreducible'' ones,   which serve as building blocks for more complicated
ones.  For example in particle physics we use irreducible representations to describe
elementary particles.   A representation \(\rho\) of a monoid on a vector space \(V\) is
\define{irreducible} if \(V\) has no subspaces invariant under all
the transformations \(\rho(m)\), except for \(\{0\}\) and \(V\) itself.  ``Irreducible 
representations'' is a bit of a mouthful, so we also call them \define{irreps} for short.

Why are irreducible representations important?   Arguably the ``indecomposable''
representations are even more important to us here.
Given two representations of a monoid, say \(\rho \maps M \to \End(V)\) and \(\rho' \maps
M \to \End(V')\), there is a representation on \(V\oplus V'\) called their \define{direct sum}:
\[     \rho \oplus \rho' \maps M \to \End(V \oplus V')  \]
given by
\[     (\rho \oplus \rho')(m)(v,v') = (\rho(m) v, \rho'(m) v') .\]
A representation is \define{indecomposable} if it is not isomorphic to a direct sum of
representations except for the 0-dimensional representation and itself.   Using an 
inductive argument we can show that every representation is a direct sum of
indecomposable representations.  That is, we can break apart any representation into
smaller pieces until we reach pieces that can't be broken apart any further.

It is easy to see that any irreducible representation is indecomposable.   The converse is not 
always true.  However, for all the monoids we shall consider here, and the kinds of representations
we consider here, indecomposability is equivalent to irreducibility!   And since ``irrep'' is such
a handy word, we shall talk about irreducibility rather than indecomposability.

\subsubsection*{\boldmath{\(S_n\)}}

Amazingly, Young diagrams can be used to classify the irreps, or at least the
``nice'' ones, of all five classical groups I listed---\(\GL(N,\C),\) \(\SL(N,\C),\) \(\U(N),\) \(\SU(N)\) 
and \(S_n\)---together with the classical monoid \(\End(\C^N)\).
Let me sketch how this goes.  We'll start with the symmetric groups \(S_n\), 
which are the most important of all.

Remember, I've shown how conjugacy classes of permutations in \(S_n\)
correspond to Young diagrams with \(n\) boxes. Now I want to
do the same for irreducible representations of \(S_n\).
This is cool for the following reason: for any finite group, the number
of irreducible representations is the same as the number of conjugacy
classes of group elements!  But in general there's no natural one-to-one
correspondence between irreducible representations with conjugacy classes. 
The group \(S_n\) just happens to be specially nice in this way.

To get started I should tell you some stuff that work for any finite group.
Suppose \(G\) is a finite group. Then \(G\) has only finitely many
irreps, all finite-dimensional. Every finite-dimensional representation
of \(G\) is a direct sum of copies of these irreps.

To get our hands on these irreps, let \(\C[G]\) be the space of
formal linear combinations of elements of \(G\). This is called the
\define{group algebra} of \(G\), since it becomes an algebra using the
product in \(G\).   With some work, one can show that \(\C[G]\) is isomorphic 
to an algebra of block diagonal matrices. For example, \(\C[S_3]\) is
isomorphic to the algebra of matrices of this form:
\[
  \left(
    \begin{array}{cccc}
      * & 0 & 0 & 0
    \\ 0 & * & 0 & 0
    \\0 & 0 & * & *
    \\0 & 0 & * & *
    \end{array}
  \right)
\] 
where the \(*\) entries can be any complex number whatsoever. Since
matrices act on vectors by matrix multiplication, we can use this to get
a bunch of representations of \(\C[G]\), and thus of \(G\) ---
one representation for each block. And this trick gives us all the
irreps of \(G\)! For example, \(S_3\) has  two \(1\)-dimensional 
irreps, coming from the two \(1\times 1\) blocks in the above
matrix, and one \(2\)-dimensional irrep, coming from the \(2\times 2\) block.

In fact, we can actually concoct these irreps as subspaces of \(\C[G]\).
One way is to find elements of \(\C[G]\) with a single 1 on the diagonal of one 
block and zero everywhere else, like these: 
\[
  \underbrace{
    \left(
      \begin{array}{cccc}
        1&0&0&0
      \\0&0&0&0
      \\0&0&0&0
      \\0&0&0&0
      \end{array}
    \right)}_{p_1}
  \qquad
  \underbrace{
    \left(
      \begin{array}{cccc}
        0&0&0&0
      \\0&1&0&0
      \\0&0&0&0
      \\0&0&0&0
      \end{array}
    \right)}_{p_2}
  \qquad
  \underbrace{
    \left(
      \begin{array}{cccc}
        0&0&0&0
      \\0&0&0&0
      \\0&0&1&0
      \\0&0&0&0
      \end{array}
    \right)}_{p_3}
\] If we can find these guys, right multiplying by them will project down
to various subspaces of \(\C[G]\), namely
\[ \{ a p_i \mid a \in\C[G]\}.\]
And these subspaces will be irreps of \(G\), as you can check using our description
of \(\C[G]\) as an algebra of block diagonal matrices.

How do we find these guys \(p_i\) in \(\C[G]\)?   That takes work!  But for starters,
notice that:
\begin{itemize}
\item
  They are \define{idempotent}: \(p_i^2 = p_i\).
\item
  They are \define{minimal}: if \(p_i\) is the sum of two idempotents,
  one of them must be zero.
\item 
 They are \define{separated}: if \(i \ne j\) we have \(p_i a p_j = 0\) for all
 \(a \in \C[G]\).
\end{itemize}
Indeed they form a large-as-possible collection of separated minimal idempotents:
as many as the number of irreps \(G\)---or equivalently, the number of conjugacy
classes in \(G\).

To go further, we need to know more about our group \(G\). So now I'll
take \(G\) to be \(S_n\) and tell you how to get 
separated minimal idempotents.  We'll get one for each Young diagram with \(n\) boxes!
Since there's as many conjugacy classes in \(S_n\) as \(n\)-box Young diagrams, 
that will mean we've got a large-as-possible collection.

Here's how it works.   Say we have a Young diagram with \(n\) boxes, like this: 
\[
\yng(3,2,1,1)
\] 
Then we can pack it with numbers from \(1\) to \(n\) like this: 
\[
\young(123,45,6,7)
\] 
There are a bunch of permutations in \(S_n\) called \define{row
permutations} that only permute the numbers within each row of our
Young diagram. And there are a bunch called \define{column permutations} that
only permute the numbers within each column.

We can form an idempotent \(p_S\) in \(\C[S_n]\) that
symmetrizes over all row permutations.
 We get \(p_S\) by taking the sum of
all row permutations divided by the number of row permutations:
\[    p_S = \frac{1}{|R|} \sum_{\sigma \in R} \sigma \in \C[S_n]  \]
where \(R\) is the set of row permutations.
Similarly, we can form an idempotent \(p_A\) in \(\C[S_n]\) that
antisymmetrizes over all column permutations. We get \(p_A\) by taking the
sum of all \emph{even} column permutations minus the sum of all
\emph{odd} column permutations, and then dividing by the total number of
column permutations:
\[   p_A = \frac{1}{|C|} \sum_{\sigma \in C}\mathrm{sgn}(\sigma) \sigma \in \C[S_n]  \]
where \(C\) is the set of column permutations.
Now here's the cool part: up to a constant factor, \(p_A q_A\) is a minimal idempotent in
\(\C[S_n]\)!  Even better, this procedure gives exactly one minimal idempotent for each
block in the block matrix description of \(\C[S_n]\).
This isn't obvious at all---it takes real work to prove---but it's the crucial fact that 
connects \(n\)-box Young diagrams to representations of \(S_n\).

Consider \(n = 3\), for example. There are 3 Young diagrams in this
case: 
\[
\yng(3)  \qquad \yng(2,1)  \qquad \yng(1,1,1)  
\] 
so \(S_3\) has 3 irreps,
confirming something I already said.   For the long squat diagram
\[   \yng(3) \]
the column permutations are trivial, so the minimal central idempotent is just \(p\).  That is,
it just ``symmetrizes'': it's the sum of all \(3!\) permutations in \(S_3\), divided by \(3!\).   It winds 
up giving a \(1 \times 1\) block in
\[  \C[S_3] \cong
  \left(
    \begin{array}{cccc}
      * & 0 & 0 & 0
    \\ 0  & * & 0 & 0
    \\0 & 0 & * & *
    \\0 & 0 & * & *
    \end{array}
  \right)
\]
and thus a 1-dimensional representation of \(S_3\).  This is the \define{trivial
representation} where every element of \(S_3\) acts as the identity operator on \(\C\).
Every monoid has a trivial representation.   

For the tall skinny diagram
\[   \yng(1,1,1) \]
the row permutations are trivial, so the minimal idempotent is just \(q\).  That is, it just 
 ``antisymmetrizes'': it's the sum of all \(3!\) permutations times their signs, divided by \(3!\).  
This gives the other 1-dimensional representation of \(S_3\):
the \define{sign representation} where each permutation acts on \(\C\) as multiplication
by its sign.  

The remaining 3-box Young diagram
\[   \yng(2,1) \]
is a bit trickier.   It gives a minimal idempotent that does a more interesting mix of row symmetrization and column antisymmetrization.   This gives the 2-dimensional representation of \(S_3\).  

Here's a more concrete way to describe this representation.  You can think of \(S_3\) as the symmetries of an equilateral triangle.  If you draw such a triangle in the plane, centered at the origin, each symmetry of this triangle gives a linear 
transformation of \(\R^2\), or in other words a \(2 \times 2\) real matrix. But you can 
think of this as a complex \(2 \times 2\) matrix!   This trick defines a
homomorphism \(\rho \maps S_3 \to \End(\C^2)\), and this is our representation.

\subsubsection*{\boldmath \(\End(\C^N)\)}

We could go on thinking about Young diagrams and representations of the symmetric
groups \(S_n\) for a long time.  People have spent their lives on this!   But before we get too
old, let's see how Young diagrams give representations of the four other classical groups.

It's actually best to start with the full linear monoid \(\End(\C^N)\), since those four
classical groups are all contained in this.   Indeed we have monoid homomorphisms like this,
all given by inclusions:
\[
\begin{tikzpicture}[scale=1.5]
\node (SU) at (0,0) {$\SU(N)$};
\node (U) at (1.5,0) {$\U(N)$};
\node (SL) at (0,-1.3) {$\SL(N,\C)$};
\node (GL) at (1.5,-1.3) {$\GL(N,\C)$};
\node (M) at (2.5,-2) {$\End(\C^N)$};
\path[->,font=\scriptsize,>=angle 90]
(SU) edge node[above]{$$} (U)
(SU) edge node[left]{$$} (SL)
(U) edge node[right]{$$} (GL)
(SL) edge node[left]{$$} (GL)
(GL) edge node[below]{$$} (M);
\end{tikzpicture}
\]
Whenever you have a monoid homorphism \(f \maps M \to M'\) and a representation of
\(M'\), say \(\rho \maps M' \to \End(V)\), you can compose them and get a representation of
\(M\).  So, representations of \(\End(\C^N\)) give representations of all four classical
groups I listed---and this is actually how we'll get our hands on irreps of these classical groups.

So let's try to understand representations of the monoid \(\End(\C^N)\).   For starters, it has an
representation on \(\C^N\) called the \define{tautologous representation}, where each transformation acts on vectors in \(\C^N\) in the obvious way.   In other words, this representation
is the identity homomorphism \(1 \maps \End(\C^N) \to \End(\C^N)\).   This is actually an irrep.

How can we get other irreps of \(\End(\C^N)\)?   One way to get new representations from old is by tensoring them.  If we have two representations \(\rho \maps M \to \End(V)\), \(\rho' \maps M \to
\End(V')\) of any monoid, we get a new one called \(\rho \otimes \rho'\) with
\[  
\begin{array}{rcl}
 \rho \otimes \rho'  \maps M &\to& \End(V) \otimes \End(V') \cong \End(V \otimes V')  \\
                                           m & \mapsto & \rho(m) \otimes \rho(m') .
\end{array}
\]                                       
So, one thing we can do is take the tautologous representation of \(\End(\C^N)\)
and tensor it with itself  a bunch of times, say \(n\) times, getting a representation on
\[\underbrace{\C^N\otimes\C^N\otimes\ldots\otimes\C^N}_{\mbox{$n\) copies}}\]
There's no reason in the world this new representation should be
irreducible. But we can try to chop it up into irreducible bits. And the
easiest way is to look for bits that transform in nice ways
when we permute the \(n\) copies of \(\C^N\). In physics lingo,
we have a space of tensors with \(n\) indices, and we can look for
subspaces consisting of tensors that transform in specified ways when we
permute the indices. For example, there will be a subspace consisting of
``totally symmetric'' tensors that don't change at all when we permute
the indices, and a subspace of ``totally antisymmetric'' tensors that
change sign whenever we interchange two indices, and so on.

But to make the ``and so on'' precise, we need Young diagrams!   After
all, these describe all the representations of the permutation group.

Here's how it works. The space
\[(\C^N)^{\otimes n} 
= \underbrace{\C^N\otimes\C^N\otimes\cdots\otimes\C^N}_{\mbox{\(n\) copies}}\]
is not only a representation of \(\End(\C^N)\); it's also a
representation of \(\C[S_n]\), coming from permutations of the \(n\) factors.  
And the actions of these two monoids commute!   This is easy to see by a direct calculation.

Next, we have seen that each \(n\)-box Young diagram diagram \(Y\) gives a minimal idempotent in \(\C[S_n]\).  This acts as an operator on \((\C^N)^{\otimes n}\), say
\[   p_Y \maps (\C^N)^{\otimes n} \to (\C^N)^{\otimes n}. \]
The image of this operator is some subspace
\[   L  = \{ p_Y v \; \big\vert \; v \in (\C^N)^{\otimes n} \} \subseteq (\C^N)^{\otimes n} .\]
But in fact, the action of  \(\End(\C^N)\) on \((\C^N)^{\otimes n}\) preserves this
subspace \(L_Y\).  Thus, \(L_Y\) becomes a representation of \(\End(\C^N)\).  So, we have
gotten a representation of \(\End(\C^N)\) from the Young diagram \(Y\)!

To see that \(L_Y\) is preserved by the action of \(\End(\C^N)\), we use the fact
that the actions of \(\End(\C^N)\) and \(\C[S_n]\) commute.   Suppose we have a 
vector in \(L\), say \(p_Y v\).  Then for any operator \(T \in \End(\C^N)\) we have
\[   T p_Y v = p_Y T v \]
so it lies in \(L_Y\).

None of this was hard. The really cool part is that \(L_Y\) is always an
 \emph{irreducible} representation of \(\End(\C^N)\).     This is much less obvious!  
The reason, ultimately, is that the linear transformations of \((\C^N)^{\otimes n}\) that commute with all transformations coming from the representation of \(\End(\C^N)\) on this space are precisely those coming from \(\C[S_n]\).   This is half of a result called ``Schur--Weyl duality''.  And I can't resist mentiong the other half, though we don't need it here.   It says that the linear transformations of \((\C^N)^{\otimes n}\) that commute with all transformations coming from the representation of \(\C[S_n]\) on this space are precisely those coming from \(\End(\C^N)\).

As you can see, there is some serious math going on here.  In any event, each Young diagram gives an irrep of \(\End(\C^N)\).   Let's see how this works in a few examples.

If we take \(n = 3\), then \(S_3\) acts on 
\[(\C^N)^{\otimes 3} = \C^N\otimes\C^N\otimes\C^N\]
So, we get some irreps of \(\End(\C^N)\) from 3-box Young diagrams.   As we've
seen, the long squat Young diagram 
\[   \yng(3) \]
gives the minimal idempotent that just ``symmetrizes''.   So it gives an irrep
of \(\End(\C^N)\) on the space of \define{symmetric tensors of rank 3}:
\[    S^3(\C^N) = \big\langle \frac{1}{3!} \sum_{\sigma \in S_n} v_{\sigma(1)} \otimes v_{\sigma(2)}
\otimes v_{\sigma(3)} \; \big\vert \; v_1, v_2, v_3 \in \C^N \big\rangle  \]
where the angle brackets mean we take all linear combinations of vectors of this form.
Similarly, the tall skinny Young diagram
\[   \yng(1,1,1)  \]
gives the minimal idempotent that ``antisymmetrizes''.   So it gives an irrep
of \(\End(\C^N)\) on the space of \define{antisymmetric tensors of rank 3}:
\[    \Lambda^3(\C^N) = \big\langle \frac{1}{3!} \sum_{\sigma \in S_n} \mathrm{sgn}(\sigma)\,
v_{\sigma(1)} \otimes v_{\sigma(2)} \otimes v_{\sigma(3)} \; \big\vert \; v_1, v_2, v_3 \in \C^N \big\rangle.  \]
All this works the same way for any other number replacing \(3\).   The other 3-box Young
diagram 
\[      \yng(2,1) \]
is more tricky.  To get its minimal idempotent up to a constant factor, you
need to first antisymmetrize over column permutations of the numbers here:
\[    \young(12,3) \]
and then symmetrize over row permutations.  Then you apply the resulting element of \(\C[S_3]\) to all vectors \(v_1 \otimes v_2 \otimes v_3\), and take all linear combinations of what you get.  
I could write down the formulas, but you probably wouldn't enjoy it.  In math, some things are more fun to do than to watch.

When you think about this game works, you'll notice that some of irreps we get are 
a bit silly. If we have a Young diagram with more than \(N\) rows, we'll be
antisymmetrizing over a tensor product of more than \(N\) vectors in \(\C^N\), which always
gives zero.  So such Young diagrams give zero-dimensional representations 
of \(\End(\C^N)\).  We can ignore these.  Indeed, most people decree that 
zero-dimensional representations don't even count as irreducible, just as the 
number \(1\) isn't prime.   Let's do that from now on.

With this convention in place, we get an irrep of \(\End(\C^N)\) from each 
Young diagram with at most \(N\) rows. 
And they're all different: that is, distinct Young diagrams with at most \(N\) rows
give nonisomorphic representations.  

Do we get \emph{all} the irreps of \(\End(\C^N)\) from Young diagrams with at most
\(N\) rows?  No, alas.   Suppose we have a representation \(\rho\) of \(\End(\C^N)\) that arises from a Young diagram.  Say it acts on some vector space \(L\).   If we pick
a basis for \(L\), we can write each linear transformation \(\rho(x) \maps L \to L\) as a matrix, and you can check that the matrix entries of \(\rho(x)\) are polynomials in the entries
of the original matrix \(x \in \End(\C^N)\).   Thus we say \(\rho\) is a \define{polynomial
representation}---and we see that Young diagrams can only give us polynomial representations
of \(\End(\C^N)\).

Thus, as soon as you find a irrep of \(\End(\C^N)\) that's not a polynomial
representation, you'll know that you can't get \emph{all} the irreps of \(\End(\C^N)\)
from Young diagrams.  And such an irrep is not hard to find.  For example, consider the
representation
\[    \begin{array}{rccc}
    \rho \maps &\End(\C^N) &\to & \End(\C^N)  \\
                               &    T  &\mapsto & \overline{T} 
 \end{array}
 \]
 that takes the complex conjugate of each entry of an \(N\times N\) matrix.   There are 
many more.

But the next best thing is true: every polynomial irrep of \(\End(\C^N)\) comes from
a Young diagram.   In fact there is a one-to-one correspondence between these
things:
\begin{itemize}
\item polynomial irreps of \(\End(\C^n)\), up to isomorphism
\item Young diagrams with \( \le N\) rows.
\end{itemize} 
Thus, we say that Young diagrams with at most \(N\) rows \emph{classify} 
polynomial irreps of \(\End(\C^N)\).   This remarkable fact is the basic link between 
Young diagrams and representations of the classical groups.   Let's see how to use it.

\subsubsection*{\boldmath\(\GL(N,\C)\)}

Let's start with the biggest of the classical groups, the general linear group
\(\GL(N,\C)\).   Consider its inclusion in \(\End(\C^N)\):
\[    \GL(N,\C) \to \End(\C^N)   \]
Composing this with any polynomial irrep of \(\End(\C^N)\), we get a representation of 
\(\GL(N,\C)\).  In fact it is an irrep.  We don't get all the irreps of \(\GL(N,\C)\), but we get all 
the polynomial irreps: that is, those whose matrix entries are polynomials in the matrix entries of the element \(g \in \GL(N,\C)\) they depend on.   

Furthermore, since \(\GL(N,\C)\) is dense in \(\End(\C^n\) and polymomials are continuous, distinct polynomial irreps of \(\End(\C^N)\) give distinct polynomial irreps of \(\GL(N,\C)\).  Even better, every polynomial irrep arises from one of \(\End(\C^n\)).  Using these ideas and our previous results on representations of \(\End(\C^N)\), we can show that there is a one-to-one correspondence between these things:
\begin{itemize}
\item polynomial irreps of \(\GL(N,\C)\), up to isomorphism
\item Young diagrams with \( \le N\) rows.
\end{itemize}
Even better, every polynomial representation of \(\GL(N,\C)\) can be written as a direct sum of polynomial irreps.

However, there are plenty of non-polynomial irreps of \(\GL(N,\C)\): not only those coming from 
the non-polynomial irreps of \(\End(\C^N)$, but also others.  The reason is that a matrix
in \(\GL(N,\C)\) has nonzero determinant, so we can cook up representations involving 
the inverse of the determinant, which is not a polynomial.   

The 1-dimensional irrep of \(\GL(N,\C)\) sending each matrix \(g\) to \(\det(g)\), called the \define{determinant representation}.  This is a polynomial irrep, so it must come from a Young diagram.  Indeed it comes from tall skinny Young diagram with one column and \(N\) rows, e.g.
\[    \yng(1,1,1,1,1)   \]
when \(N = 5\).    If we have any irrep of \(\GL(N,\C)\) coming from a Young diagram,
tensoring it with the determinant representation gives a new irrep described by a Young
diagram with an extra column with \(N\) rows, like this:
\begin{center}
\( \yng(4,3,1,1,1) \) \raisebox{2.5em}{\(\quad \otimes \qquad\)} \(\yng(1,1,1,1,1)\)
\raisebox{2.5em}{\qquad \(\cong\) \qquad} \( \yng(5,4,2,2,2) \)
\end{center}
However, there's also a 1-dimensional irrep of \(\GL(N,\C)\) that sends \(g \in \GL(N,\C)\) to \(\det(g)^{-1}\).   This is called the \define{inverse} of the determinant representation, both for the obvious reason and because when you tensor it with the determinant representation you get the trivial representation.   Since \(\det(g)^{-1}\) is not a polynomial in the matrix entries of \(g\), this not a polynomial representation.    But it is still an \define{algebraic representation}: one whose matrix entries are rational functions of the matrix entries of \(g\).   

Algebraic representations are the kind most natural in algebraic geometry.   
Indeed \(\GL(N,\C)\) is a \define{linear algebraic group} over \(\C\): that is,
a group in the category of affine algebraic varieties over the complex numbers.   When people talk about representations of linear algebraic groups, they usually mean algebraic representations.   

So, fans of algebraic geometry will be glad to know that algebraic irreps of \(\GL(N,\C)\) 
can all be built by taking a polynomial irrep and tensoring it with the inverse of the 
determinant representation some number of times.   This in turn means we can describe any algebraic irrep of \(\GL(N,\C)\) using a Young diagram with \emph{fewer than} \(N\) rows 
together with an integer \(k\).   The Young diagram gives a representation \(\rho\), and 
then we form the representation on the same space where \(g\) acts by \( \det(g)^k \rho(g) \).   
If \(k \ge 0\) this is the same as tacking on \(k\) extra columns with \(N\) rows to our 
Young diagram, but the procedure also makes sense for \(k < 0\).   We get a one-to-one
correspondence between these things:
\begin{itemize}
\item algebraic irreps of \(\GL(N,\C)\), up to isomorphism
\item pairs consisting of a Young diagram  with \(< N\) rows and an integer.
\end{itemize}
If you like, you can think of such a pair as a funny sort of Young diagram with \(\le N\)
rows where the number of columns with \(N\) rows can be any integer---even a negative number!

This is the story for irreps, but what about more general representations?  It's as nice as
it could be: every algebraic representation of \(\GL(N,\C)\) is a direct sum of algebraic
irreps.

If you don't yet love algebraic geometry, you may prefer to think of \(\GL(N,\C)\) as a 
\define{complex Lie group}: a group in the category of complex manifolds.    When we talk about a representation of a complex Lie group \(G\), we usually mean an \define{complex-analytic representation}: a representation
\(\rho \maps \GL(N,\C) \to \End(L)\) for which the matrix entries of
\(\rho(g)\) are complex-analytic functions of the matrix entries of \(g\).   
Luckily for \(\GL(N,\C)\)  these representations are all algebraic!  The constraint 
\(\rho(gh) = \rho(g)\rho(h)\) is so powerful that any complex-analytic solution is
actually algebraic.  So, the whole story we told for algebraic representations of
\(\GL(N,\C)\) also applies to complex-analytic ones.

\subsubsection*{\boldmath\(\SL(N,\C)\)}

We can also get representations of the special linear group \(\SL(N,\C)\) from Young diagrams.  Any Young diagram with at most \(N\) rows gives an algebraic irrep of \(\End(\C^N)\), and composing this with the inclusion
\[    \SL(N,\C) \to \End(\C^N)   \]
we get an algebraic irrep of \(\SL(N,\C)\).   We get all the algebraic irreps of \(\SL(N,\C)\)
this way.   Even better, the irritating fly in the ointment for \(\GL(N,\C)\), the determinant representation, become trivial for \(\SL(N,\C)\).   So does the inverse of the determinant
representation.   So, we get a one-to-one correspondence between these two things:
\begin{itemize}
\item algebraic irreps of \(\SL(N,\C)\), up to isomorphism
\item Young diagrams with \( < N\) rows.
\end{itemize}
Furthermore, every algebraic representation of \(\SL(N,\C)\) is a direct sum of algebraic irreps.
So, algebraic representations of \(\SL(N,\C)\) are classified by \emph{finite collections} of
Young diagrams with \( < N\) rows.

Here we are thinking of \(\SL(N,\C)\) as a linear algebraic group.    We can also think of it as a complex Lie group.   However, all its complex-analytic representations are algebraic.  So the same classification applies here too.  

\subsubsection*{\boldmath\(\U(N)\)}

The unitary group \(\U(N)\) is different from the classical groups so far, because the equations defining unitarity involve complex conjugation:
\[         g g^\ast = 1 \]
so it's not a linear algebraic group over \(\C\).   Instead it's a linear algebraic group over \(\R\).   We  shall still study its representations on complex vector spaces, but now the interesting ones are
the \define{real-algebraic representations}:  those where the matrix entries of \(\rho(g)\) are 
rational functions of the real and imaginary parts of the matrix entries of \(g\).

To get representations of \(\U(N)\) it's convenient to use our knowledge of representations of
\(\GL(N,\C)\).    We can take any algebraic irrep of \(\GL(N,\C)\) and compose it with the inclusion
\[    \U(N) \to \GL(N,\C)  \]
to get a real-algebraic representation of \(\U(N)\).   The result is an irrep, and we get all the 
real-algebraic irreps of \(\U(N)\) on complex vector spaces this way.   In fact, the classification of these real-algebraic irreps of \(\U(N)\) completely matches the classification of algebraic irreps of \(\GL(N,\C)\).     We thus get a one-to-one correspondence between these things:
\begin{itemize}
\item real-algebraic irreps of \(\U(N)\) on complex vector spaces, up to isomorphism
\item pairs consisting of a Young diagram with \(< N\) rows and an integer.
\end{itemize}
Furthermore, every real-algebraic representation of \(\U(N)\) is a direct sum of real-algebraic
irreps.

Alternatively, we can think of \(\U(N)\) as a \define{Lie group}: a group in the category of
manifolds (ordinary real manifolds, not complex manifolds).   For a Lie group it's natural to study
\define{smooth representations}:  those where the matrix entries of \(\rho(g)\) are smooth functions of the matrix entries of \(g\).   Or we can go further and think of \(\U(N)\) as a mere
\define{topological group}: a group in the category of topological spaces.   For a topological
group it's natural to study \define{continuous representations}, where the matrix entries of \(\rho(g)\) are continuous functions of the matrix entries of \(g\).  

But something very nice is true: every smooth representation of \(\U(N)\) is automatically real-algebraic, and every continuous representation of \emph{any} Lie group is automatically smooth!   So we do not gain any generality by considering smooth or continuous irreps of \(\U(N)\): they are both classified by pairs consisting of a Young diagram with \(< N\) rows and an integer.

Another variant also turns out to work the same way.  In quantum physics we use unitary representations on Hilbert spaces.  A finite-dimensional Hilbert space, which is the only kind we'll consider here, is just a finite-dimensional complex vector space with an inner product.  
A \define{unitary representation} of a group \(G\) on a Hilbert space \(H\) is a  representation \(\rho \maps G \to \End(V)\) such that each of the transformations \(\rho(g)\) is unitary.

 It turns out that because \(\U(N)\) is compact, we can take any continuous representation \(\rho \maps \U(N) \to \End(V)\), pick any inner product on the vector space \(V\), and ``average it'' over the action of \(U(N)\) to get a new improved inner product with
\[          \langle \rho(g) v, \rho(g) w \rangle = \langle v, w \rangle  \quad \textrm{ for all }
v, w \in V \textrm{ and } g \in \U(N) .\]
This says that all the transformations \(\rho(g)\) are unitary:
\[        \rho(g)^* \rho(g) = 1 .\]
So, \(\rho\) has been promoted to a unitary representation.   

Putting this together with what we already have, one can show there is
a one-to-one correspondence between these things:
\begin{itemize}
\item continuous unitary irreps of \(\U(N)\), up to isomorphism
\item pairs consisting of a Young diagram with \(< N\) rows and an integer.
\end{itemize}
Also, every continuous unitary representation of \(\U(N)\) is a direct sum of continuous unitary irreps.

\subsubsection*{\boldmath\(\SU(N)\)}

Finally we turn to the special unitary group \(\SU(N)\).  Since all the main patterns have been laid out, we will go faster now---as usual, not proving things but at least trying to make them plausible.  Just as \(GL(N,\C)\) helps us understand \(\U(N)\), \(\SL(N,\C)\) helps us understand \(\SU(N)\).    The reason, ultimately, is that \(\U(N)\) is  the``compact real form'' of the complex Lie group  \(\GL(N,\C)\), and \(\SU(N)\) is the compact real form of \(\SL(N,\C)\).   But to understand this, one needs to get into Lie theory more deeply than we intend to here.

We can take any algebraic irrep of \(\SL(N,\C)\) and compose it with the inclusion
\[    \SU(N) \to \GL(N,\C)  \]
to get a representation of \(\SU(N)\).  This is a real-algebraic irrep, and we get all the 
real-algebraic irreps of \(\SU(N)\) this way.   With help from our classification of algebraic irreps of  \(\SL(N,\C)\), we we can show there is a one-to-one correspondence between these things:
\begin{itemize}
\item real-algebraic irreps of \(\SU(N)\), up to isomorphism
\item Young diagrams with \(< N\) rows.
\end{itemize}
Then, by the averaging trick mentioned already for \(\U(N)\), we also get a one-to-one correspondence between these things:
\begin{itemize}
\item continuous unitary irreps of \(\SU(N)\), up to isomorphism
\item Young diagrams with \(< N\) rows.
\end{itemize}
Further more, as we have come to expect, in both the real-algebraic case and the continuous unitary case every representation of the given sort is a direct sum of irreps of that sort.

\subsubsection*{Summary and further directions}

Let's summarize what we have seen---but also say a bit more.  While we have studied representations on finite-dimensional vector spaces over \(\C\), most of the \emph{purely
algebraic} results hold for any field of characteristic zero!   Fields with nonzero characteristic
behave very differently, and in fact the irreducible representations of \(S_n\) still haven't
been classified over finite fields.  But the items with check marks here hold if we replace \(\C\) with any field of characteristic zero:

\begin{itemize}
\item[\ding{51}] Irreps of \(S_n\) are classified by Young diagrams with \(n\) boxes.
\item[\ding{51}] Polynomial irreps of \(\End(\C^N)\) are classified by Young diagrams with 
\(\le N\) rows.
\item[\ding{51}] Polynomial irreps of \(\GL(N,\C)\) are classified by Young diagrams with
\( \le N\) rows.  
\item[\ding{51}] Algebraic irreps of \(\GL(N,\C)\) are classified by pairs consisting of a Young diagram with \( < N\) rows and an integer.
\item[\ding{51}] Algebraic irreps of \(\SL(N,\C)\) are classified by Young diagrams with \(< N\) rows.
\item Analytic irreps of \(\SL(N,\C)\) are classified by Young diagrams with \(< N\) rows.
\item Analytic irreps of \(\GL(N,\C)\) are classified by pairs consisting of a Young diagram 
with \( < N\) rows and an integer.
\item Real-algebraic irreps of \(\U(N)\) are classified by pairs 
consisting of a Young diagram with \(< N\) rows and an integer. 
\item Continuous unitary irreps of \(\U(N)\) are classified by pairs 
consisting of a Young diagram with \(< N\) rows and an integer.  
\item Real-algebraic irreps of \(\SU(N)\) are classified by 
Young diagram with \(< N\) rows.
\item Continuous unitary irreps of \(\SU(N)\) are classified by 
Young diagram with \(< N\) rows.
\end{itemize}

However, this is far from the end of the story!
First of all, we can use $n$-box Young diagrams packed with numbers \(1,\dots,n\), called
\define{Young tableaux}, to do all sorts of calculations involving irreps of classical
groups.   

Say we want to figure out the dimension of the irrep of \(S_n\) corresponding
to some Young diagram.  Then we just count the \define{standard Young tableaux}
of that shape: that is, Young tableaux where the numbers increase as we go down any
column or across any row.   For example, there are two standard Young tableaux
of this shape:
\[   \young(12,3) \qquad \young(13,2) \]
so this Young diagram:
\[    \yng(2,1) \]
gives a 2-dimensional irrep of \(S_3\).

Or: say we tensor two irreps and want to
decompose the result as a direct sum of irreps: how do we do it? We play a
little game with Young tableaux and out pops the answer.  The relevant
buzzword is ``Littlewood--Richardson rules''. Or say we have an irrep of
\(S_n\) and want to know how it decomposes into irreps when we restrict
it to a subgroup like \(S_{n-1}\), or similarly for \(\SL(N,\C)\)
and \(\SL(N-1,\C)\), etc. How do we do this? More messing with Young
tableaux. Here one relevant buzzword is ``branching rules''.

I'll warn you right now: there is an \emph{enormous} literature on this
stuff. The combinatorics of Young diagrams is one of those things that
everyone has worked on, from hardnosed chemists to starry-eyed category
theorists. It takes a lifetime to master this material, and I certainly
have \emph{not}. But learning even a little is fun, so don't be
\emph{too} scared.

Second of all, Young diagrams are also good for studying the
representations of some other classical groups, such as these:
\begin{itemize}
\item The \define{orthogonal group} \(\mathrm{O}(N)\), consisting of all 
orthogonal linear transformations of \(\R^N\).
\item The \define{special orthogonal group} \(\mathrm{SO}(N)\), consisting
of all orthogonal linear transformations of \(\R^N\) with determinant \(1\).
\item The \define{symplectic group} \(\mathrm{Sp}(2N)\), consisting of all
symplectic linear transformations of \(\R^{2N}\).
\end{itemize}
All these groups have an
obvious ``tautologous representation'', and we can cook up other representations
by taking the \(n\)th tensor power of this representation and
hitting it with minimal idempotents in \(\C[S_n]\) coming from Young diagrams. The
story I just told you can be repeated with slight
or not-so-slight variations for these other groups.

Third, we can ``\(q\)-deform'' the whole story, replacing any one of
these classical groups by the associated ``quantum group'', and
replacing \(\C[S_n]\) by the corresponding ``Hecke algebra''.
This is really important in topological quantum field theory and the
theory of von Neumann algebras.

Fourth, there are nice relationships between Young diagrams and
algebraic geometry, like the ``Schubert calculus'' for the cohomology
ring of a Grassmannian.

Fifth and finally, Young diagrams are themselves objects in an important 
category!

To understand this we need to step back a bit.   We have seen that Young diagrams 
are good for getting new representations from old ones.  Given any
representation 
\[  \rho \maps M \to \End(V) \]
of any monoid \(M\), and given any Young diagram \(Y\), we can get a new 
representation of \(M\) as follows.  First form the \(n\)th tensor power of \(\rho\),
which is the representation
\[   \rho^{\otimes n} \maps M \to \End(V^{\otimes n}) \]
defined by
\[   \rho^{\otimes n}(m)(v_1 \otimes \cdots \otimes v_n) = 
\rho(m)(v_1) \otimes \cdots \otimes \rho(m)(v_n). \]
The group \(S_n\) also acts on \(V^{\otimes n}\), so the minimal idempotent in 
\(\C[S_n]\) coming from \(Y\) gives an idempotent operator
\[       p_Y \maps V^{\otimes n} \to V^{\otimes n} \]
Then take the image of $p_Y$.   Since the actions of $M$ and $S_n$ on $V^{\otimes n}$
commute, this image is a subspace of \(V^{\otimes n}\) that is invariant under all the transformations \(\rho(m)\) for \(m \in M\).   So, it gives a representation of \(M\).
Let us call this new representation \(Y(\rho)\).

Since this procedure for getting new representations from old is completely systematic, it should be a functor. 
Indeed, this is true!  There is a category $\Rep(M)$ whose objects are representations of $M$, with the usual morphisms between these.   There is a
functor from this category to itself, say
\[    Y \maps \Rep(M) \to \Rep(M) , \]
that maps each representation $\rho$ to $Y(\rho)$.   And this functor is
called a \define{Schur functor}.

Schur functors also work on categories other than categories of representations.   
Very roughly, Schur functors know how to act on any category where:
\begin{itemize}
\item we can take linear combinations of morphisms $f, g \maps x \to y$
between any two objects $x$ and $y$,
\item we can take direct sums and tensor products of objects,
\item the symmetric group $S_n$ acts on $x^{\otimes n}$ for any object $x$, and
\item we can project to the image of any idempotent morphism $f \maps x \to x$.
\end{itemize}
One can make these conditions precise, and I have taken to calling categories obeying these conditions ``2-rigs''.   So, for any 2-rig  $\RR$ and any Young diagram $Y$, we get a Schur functor
\[      Y_\RR \maps \RR \to \RR . \]
(Now I am being more careful to indicate that the Schur functor depends on the category \(\RR\).)

There is a nice way think about what is going 
on here.  There is a 2-rig $\Schur$ whose objects are formal finite direct sums of Young
diagrams, like this:
\[
 \Yvcentermath1  \yng(1,1,1,1) \;\; \oplus \;\; \yng(4,2,1) \;\; \oplus \;\; \yng(2,2) \;\; \oplus \;\; \yng(1) \;\;\oplus \;\; \yng(1)
\]
This 2-rig \(\Schur\) plays a special role in the theory of 2-rigs: it is the ``free 2-rig on one object''.   This object is the one-box Young diagram:
\[    \yng(1) \]
What does this mean?  Roughly speaking, it means that for any 2-rig \(\RR\) and any object \(r \in \RR\), there is a unique functor (or more precisely, map of 2-rigs)
\[       F \maps \Schur \to \RR \]
sending the one-box Young diagram to \(r\):
\[     \Yvcentermath1    F(\,\yng(1)\,)  = r .\] 
This functor \(F\) must send each Young diagram \(Y\) to some object in \(\RR\).  Which object is that?   It is the result of applying the Schur functor corresponding to \(Y\) to \(r\):
\[       F(Y)  = Y_\RR(r) . \]
While these ideas may seem painfully abstract, they are elegant, and they turn out to clarify many topics in the theory of Young diagrams---see the references for more details.

\subsubsection*{References}

I have zipped through a lot of material but not explained it in detail.   The lectures I gave at the University of Edinburgh, based on these notes, may help:
\begin{itemize}
\item 
John C.\ Baez, \href{https://www.youtube.com/playlist?list=PLuAO-1XXEh0a4UCA-iOqPilVmiqyXTkdJ}{Talks on This Week's Finds
in Mathematical Physics}, Lectures 1--3.
\end{itemize}
But there are still many details missing.  So, how can you really learn this stuff? 

If you have a certain amount of patience for old-fashioned terminology, I recommend going back to the classic text on classical groups:
\begin{itemize}
\item
  Hermann Weyl, \emph{The Classical Groups, Their Invariants and
  Representations}, Princeton U. Press, Princeton, 1997.
\end{itemize}
Weyl coined the term ``classical groups'' for the purposes of this book,
which was first published in 1939. His prose is beautiful, but I warn
you, this book is not the way to learn Young diagrams in a hurry.

For a user-friendly approach that's aimed at physicists, but still
includes proofs of all the key results, you can't beat this:
\begin{itemize}
\item
  Irene Verona Schensted, \emph{A Course on the Applications of Group
  Theory to Quantum Mechanics}, NEO Press, Box 32, Peaks Island, Maine.
\end{itemize}
A girlfriend gave me a copy when I was a college student, but
only much later did I realize how great a book it is. Unfortunately it's
out of print! Someone should reprint this gem.   In the meantime, here is 
another book that covers Young diagrams and their applications to physics:
\begin{itemize}
\item
  Shlomo Sternberg, \emph{Group Theory and Physics}, Cambridge U. Press,
  Cambridge, 1994.
\end{itemize}
Both these books, but especially the latter, describe applications of
Young diagrams to particle physics, like Gell-Mann's famous ``eight-fold
way'', which was based on positing an \(\SU(3)\) symmetry
between the up, down and strange quarks.

Then there are more advanced texts, for when your addiction to Young
diagrams becomes more severe. For the combinatorial side of things,
these are good:
\begin{itemize}
\item
  Gordon Douglas James and Adalbert Kerber, \emph{The Representation
  Theory of the Symmetric Group}, Addison-Wesley, Reading,
  Massachusetts, 1981.
\item
  Bruce Eli Sagan, \emph{The Symmetric Group: Representations,
  Combinatorial Algorithms, and Symmetric Functions}, Springer, Berlin,
  2001.
\end{itemize}
For a more conceptual approach to representation theory that puts Young
diagrams in a bigger context, try these:
\begin{itemize}
\item William Fulton and Joe Harris, \emph{Representation Theory --- a First Course}, Springer, Berlin, 1991. 
\item
  Roe Goodman and Nolan R.\ Wallach, \emph{Representations and Invariants
  of the Classical Groups}, Cambridge University Press, Cambridge, 1998.
\end{itemize}
And finally, here's a
mathematically sophisticated book that really gives you a Young diagram
workout:
\begin{itemize}
\item
  William Fulton, \emph{Young Tableaux: With Applications to
  Representation Theory and Geometry}, Cambridge University Press, Cambridge,
  1997.
\end{itemize}

If you want to learn about Lie groups, there are many good books.  I'll list some in rough order
of increasing sophistication:
\begin{itemize}
\item Brian Hall, \emph{Lie Groups, Lie Algebras, and Representations}, Springer, Berlin, 2003.
\item J.\ Frank Adams, \emph{Lectures on Lie Groups}, University of Chicago Press, Chicago, 2004.
\item Sigurdur Helgason, \emph{Differential Geometry, Lie Groups, and 
Symmetric Spaces}, Academic Press, New York, 1979.
\item  Daniel Bump, \emph{Lie Groups}, Springer, Berlin, 2004. 
\end{itemize}
The book by Fulton and Harris starts with an introduction to representations of 
finite groups, especially \(S_n\), and it has a chapter on Young diagrams.  For linear algebraic groups, try this:
\begin{itemize}
\item
James S.\ Milne, \emph{\href{https://math.ucr.edu/home/baez/qg-fall2016/Milne_AGS.pdf}{Basic Theory of Affine Group Schemes}}.
\end{itemize}

Finally, this paper explains how the category \(\Schur\), whose objects are formal finite direct sums of Young diagrams, is the free 2-rig on one object:
\begin{itemize}
\item
  John C.\ Baez, Joe Moeller and Todd Trimble, \href{https://arxiv.org/abs/2106.00190}{Schur       
  functors and categorified plethysm}.
\end{itemize}
There is a known way to compose formal direct sums of Young diagrams, called ``plethysm'', and we study plethysm using the 2-rig \(\Schur\).

\subsection*{Acknowledgements}

I thank the Leverhulme Trust for giving me a fellowship to give a series of lectures on this topic at the School of Mathematics of the University of Edinburgh, and I thank Tom Leinster for making this actually happen.  I thank James Dolan, Joe Moeller and Todd Trimble for many discussions of Young diagrams.

\end{document}